\documentclass[preprint,12pt]{elsarticle}
\usepackage{amsmath,amsfonts}
\usepackage{algorithmic}
\usepackage{algorithm}
\usepackage{array}
\usepackage[caption=false,font=normalsize,labelfont=sf,textfont=sf]{subfig}
\usepackage{textcomp}
\usepackage{stfloats}
\usepackage{url}
\usepackage{verbatim}
\usepackage{graphicx}
\usepackage{amsmath,amssymb,amsfonts}
\usepackage{algorithmic}
\usepackage{graphicx}
\usepackage{textcomp}
\usepackage{bm}

\usepackage{amsthm}

\theoremstyle{definition}

\theoremstyle{proposition}

\newcommand{\norm}[1]{\left\lVert#1\right\rVert}

\begin{document}
	
\begin{frontmatter}
	
	\title{Optimization of Fossil Fuel Consumption under Grid Energy Supply by Wind Energy and Auxiliary Batteries Using Fossil Fuel for Energy Supply}
	
	\tnotetext[t1]{This paragraph of the first footnote will contain the date on 
		which you submitted your paper for review. }
	
	\author[1]{Yi Chen \corref{cor1}}
	\ead{chenyi@hust.edu.cn}

	\address[1]{The China-EU Institute for Clean and Renewable Energy, Huazhong University of Science and Technology, Wuhan, Hubei, China}
	\cortext[cor1]{Corresponding author}

	\begin{abstract}
		In this paper, a class of systems in which auxiliary batteries fed by fossil fuel generation and wind farms work together to supply power to the grid are modeled and strategies are solved to minimize the consumption of fossil fuels while still allowing the total energy to satisfy the energy demand of the grid. The wind resource is modeled using stochastic differential equations and the grid energy demand is also stochastically modeled using stochastic differential equations. The two stochastic models are well calibrated using wind resource data from historical wind farms and historical grid energy demand data. Finally, dynamic programming is used to solve this optimization problem.
	\end{abstract}
	
\end{frontmatter}

	\section{INTRODUCTION}
	The issue of energy has always been a huge problem for the survival of mankind. From the earliest period when mankind used biomass as the main source of energy, to the time when mankind could use electricity, the innovations have been tremendous \cite{smil2004world}. However, the current level of technology still does not allow humans to consume energy at will. In addition, with the progress of time, in spite of the problem of the total amount of energy, the serious impact on the environment caused by the massive use of fossil fuels has also aroused the concern of mankind. Humans began to try to find a larger total amount of cleaner energy.
	Renewable energy sources have received a great deal of attention in recent years due to mankind's keen interest in cleaner, more extensive (and in some sense inexhaustible) energy sources. Among them, solar and wind energy are the most favored children of renewable energy\cite{vakulchuk2020renewable}.
	Solar and wind energy are two of the most widely available sources of renewable energy and are the most accessible to mankind at the current state of the technology. For these two types of energy, there are already a large number of researchers working in their fields, such as solar thermal power generation and solar photovoltaic power generation using solar energy. Wind turbine research utilizing wind energy\cite{qazi2019towards}.
	The use of both types of energy is clean and in some sense inexhaustible. However, due to the two energy sources have unstable factors in the utilization process, such as solar energy can not be directly used for power generation at night\cite{kuang2016review}, wind energy can not guarantee a stable power generation at any time\cite{gross2003progress}. Therefore, there is still a lot of work to be done by researchers to try to utilize new energy instead of conventional energy sources.
	
	There is a great deal of research on the uncertainty of wind energy, such as\cite{kwon2010uncertainty}\cite{moehrlen2004uncertainty}\cite{lackner2007uncertainty}\cite{yan2019uncertainty}. Modeling of uncertainty in wind energy has been done using neural networks for wind speed prediction\cite{zhou2019wind}\cite{banik2020uncertain}\cite{zhang2019short} and stochastic differential equations for wind speed assessment\cite{moller2016probabilistic}\cite{verdejo2019modelling} \cite{iversen2016short}. The energy demand of the grid has been studied in\cite{almansoori2012design}\cite{de2022predictive}\cite{yang2021robust}
	
	The structure of this paper is modeled using stochastic differential equations for wind resources and grid energy demand in \ref{problem}, where a common calibration method is provided for stochastic differential equations by means of a large amount of historical data in \ref{calibration}. In \ref{optimal}, an optimization problem is created for this established stochastic model. A common method for solving this optimization problem is given in \ref{optimal}. Simulation results are given in \ref{simulation}.
	\section{PROBLEM MODELING FOR WIND ENERGY AND ENERGY DEMADN}
	\label{problem}
	For a fixed area, the presence or absence of wind and the magnitude of wind speed are influenced by a large number of uncertainties in nature. With the current level of human technology, it is not possible to establish a complete mathematical model to give the size of the wind speed at any time in a certain area, which has a great impact on the utilization of wind energy.
	For this reason, this paper will use stochastic differential equations to establish a stochastic model for the wind speed in a certain area, as shown in (\ref{wind speed sde})
	\begin{equation}
		\label{wind speed sde}
		\dot{v^w_t} = f^{\theta_w}(v^w_t)dt+g^{\theta_w}(v^w_t)dB^w
	\end{equation}
	where $v^w_t$ is the  wind speed at time $t$ in a certain area. $f^{\theta_w}$ and $g^{\theta_w}$ represent the drift and diffusion terms in the stochastic differential equation, respectively, and both functions will be fitted using a neural network through past wind resource data for the region. $B^w$ is a $1$-dimensional Brownian motion.
	
	For wind energy at moment $t$ can be expressed by (\ref{wind energy})
	\begin{equation}
		\label{wind energy}
		P^w_t = \frac{1}{2}\rho A_w (v^w_t)^3
	\end{equation}
	where $\rho$ and $A_w$ are the density of air and the area through which the wind passes at that wind speed, respectively.
	According to It\^{o}'s chain rule (\cite{evans2012introduction}section4.4), the following equation holds
	\begin{equation}
		\begin{aligned}
			\mathbb{E}&\left\{\frac{dP^w_t}{dt}\right\} = 
			&\frac{3}{2} \rho A_w \left(v^w_t\right)^2 f^{\theta_w}\left(v^w_t\right)+ 3\rho A_w v^w_t \left(g^{\theta_w}\left(v^w_t\right)\right)^2 
		\end{aligned}
	\end{equation}
	The amount of energy that can be obtained from wind energy over a period of time can be described by the following equation.
	\begin{equation}
		\begin{aligned}
			E^{t_0,T}_w = \int_{t_0}^{T}\Bigg(\int_{t}^{t+dt}\bigg(\frac{3}{2} \rho A_w &\left(v^w_s\right)^2 f^{\theta_w}\left(v^w_s\right)+\\
			 &3\rho A_w v^w_s \left(g^{\theta_w}\left(v^w_s\right)\right)^2 \bigg)ds\Bigg)dt
		\end{aligned}
	\end{equation}
	where $t_0$ is the initial time for the observation of wind speed and $T$ is the termination time.
	Likewise, there is a large amount of uncertainty in the energy supplied by the grid to the population, since the electricity demand of the population is not a continuous constant in the time domain. Therefore, a stochastic model can be developed for the energy consumption demand of the residents using stochastic differential equations as well.
	\begin{equation}
		\dot{E^d_{t}} = \mu^{\theta_d}(E^d_{t})dt+\sigma^{\theta_d}(E^d_{t})dB^d
	\end{equation}
	where $\mu^{\theta_d}$ and $\sigma^{\theta_d}$ represent the drift and diffusion terms in the stochastic differential equation, respectively, and both functions will be fitted using a neural network through past energy demand data for the region. $B^d$ is a $1$-dimensional Brownian motion.
	\subsection{Calibration of Stochastic Differential Equations}
	\label{calibration}
	In this section, we use a more common method to calibrate the stochastic differential equation for wind speed and the stochastic differential equation for energy demand, respectively.
	\begin{equation}
		\begin{aligned}
			\mathcal{L}(\theta_w) = \int_{t_0}^{T} p(v^w_t\vert \theta_w)dt,\\
			\mathcal{L}(\theta_d) = \int_{t_0}^{T} p(E^d_t\vert \theta_d)dt
		\end{aligned}
	\end{equation}
	where $p(v^w_t\vert \theta_w)$ and $p(E^d_t\vert \theta_d)$ are the probability of occurrence of $v^w_t$ under $theta_w$ parameter and the probability of occurrence of $E^d_t$ under $theta_d$ parameter, respectively. The purpose of the calibration is to maximize $\mathcal{L}(\theta_w)$ and $\mathcal{L}(\theta_d)$ by adjusting the parameters $theta_w$ and $theta_d$.
	In addition, the stochastic differential equation has a solution needs to ensure that its drift term function and diffusion term function satisfy the Lipschitz continuous condition, so the total loss functions of $\mathcal{L}(\theta_w$) and $\mathcal{L}(\theta_w)$ are the following equations, respectively.
		\begin{equation}
		\begin{aligned}
			\mathcal{L}(\theta_w) &= -\int_{t_0}^{T} \log p(v^w_t\vert \theta_w)dt\\
			&+ \kappa_1 \int_{t_0}^{T} \left(\norm{\partial f^{\theta_w}(v^w_t)}_2 + \norm{\partial g^{\theta_w}(v^w_t)}_2-C_1\right)^2dt,\\
			\mathcal{L}(\theta_d) &= -\int_{t_0}^{T} \log p(E^d_t\vert \theta_d)dt\\
			&+ \kappa_2 \int_{t_0}^{T} \left(\norm{\partial \mu^{\theta_d}(E^d_t)}_2 + \norm{\partial \sigma^{\theta_d}(E^d_t)}_2-C_2\right)^2dt,
		\end{aligned}
	\end{equation}
	where $\kappa_1$ and $\kappa_2$ are positive hyperparameters that set the weight of the regular term in the loss function for Lipschitz continuous condition. $C_1$ and $C_2$ are arbitrary hyperparameters.
	\section{Auxiliary Battery Optimization Energy Consumption Problem Establishment and Solving}
	\label{optimal}
\begin{figure}[htbp]
		\centering
		\includegraphics[width=\columnwidth]{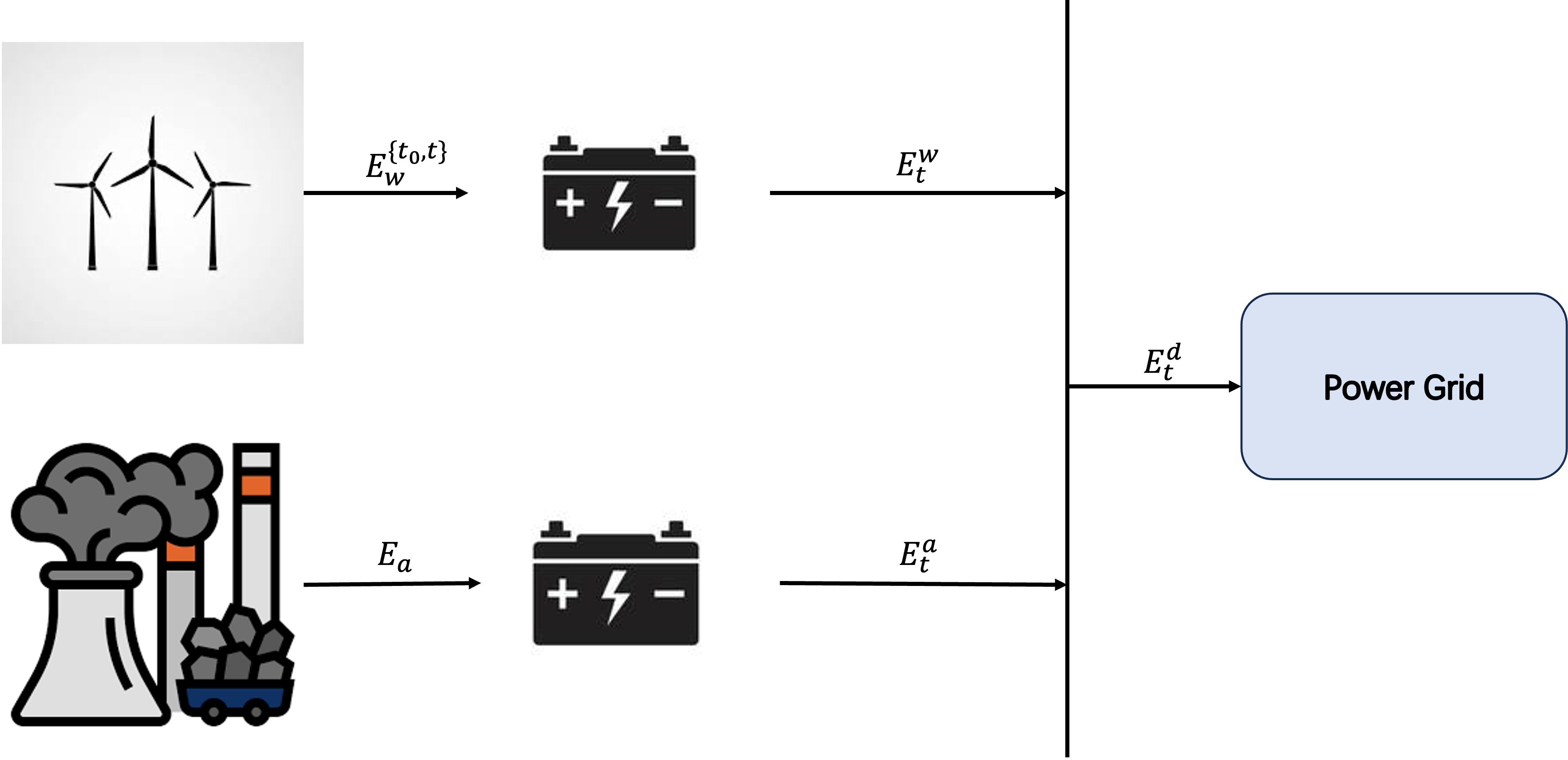}
		\caption{Wind power and auxiliary batteries for grid energy supply}
		\label{energy optimal}
	\end{figure}

	The energy of the auxiliary battery serves as the energy used to supplement the energy generated by the wind power when it is not sufficient to supply the energy demand of the grid and the source of this energy is generated from fossil energy generation.
	The goal of this section is to minimize the use of auxiliary battery energy, provided that wind and auxiliary batteries can supply the energy needs of the grid.
	\begin{equation}
		\begin{aligned}
			&\min E^a_t\\
			s.t.\quad &E^d_t = E^w_t + E^a_t
		\end{aligned}
	\end{equation}
	where $E^a_t$ is the energy that can be supplied by the auxiliary battery at moment $t$, and $E^w_t$ is the energy that can be supplied by the wind energy stored in the battery at moment $t$, as shown in Fig. \ref{energy optimal}.
	For simplicity of the model, $E^w_t$ and $E^a_t$ can be considered to obey the following equation.
	\begin{equation}
		\begin{aligned}
			E^w_t = \eta_w E^{t_0,t},\\
			E^a_t = \eta_a E_a
		\end{aligned}
	\end{equation}
	where $E^{t_0,t}$ can be considered as the total energy produced by wind energy until moment $t$ and $E_a$ can be considered as the total energy produced by fossil fuels.
	The optimization problem can be solved by treating $E^a_t$ as a control input for this optimization problem.

	\section{SIMULATION}
	\label{simulation}

	\bibliography{WEM.bib}
	\bibliographystyle{elsarticle-num}
	
\end{document}